\documentclass[12pt]{article}

\usepackage{a4wide}

\usepackage[intlimits]{amsmath}
\usepackage{amssymb,amsthm}

\usepackage[cp1251]{inputenc}
\usepackage[english]{babel}

\newtheorem{thm}{Theorem}[section]
\newtheorem{cor}[thm]{Corollary}
\newtheorem{lem}[thm]{Lemma}

\theoremstyle{definition}

\theoremstyle{remark}
\newtheorem{rem}[thm]{Remark}
\numberwithin{equation}{section}

\begin{document}\Large

\title{Asymptotic expansion of semi-Markov random evolutions}

\author{S. Albeverio$^{1,2,3}$, V.S. Koroliuk$^{4}$, I.V.Samoilenko$^{4,*}$\\ $^1$Institut f\"{u}r Angewandte Mathematik,
Universit\"{a}t Bonn,\\ Wegelerstr. 6, D-53115 Bonn
(Germany)\\$^2$SFB 611, Bonn, BiBoS, Bielefeld-Bonn\\ $^3$CERFIM,
Locarno and USI (Switzerland)\\ $^4$Institute of Mathematics,
Ukrainian National Academy of
Sciences,\\ 3 Tereshchenkivs'ka, Kyiv, 01601, Ukraine\\
$^{*}$isamoil@imath.kiev.ua (corresponding author)}

\maketitle

\abstract{Regular and singular parts of asymptotic expansions of
semi-Markov random evolutions are given. Regularity of boundary
conditions is shown. An algorithm for calculation of initial
conditions is proposed. 
\vskip 3mm {\bf AMS-Classification:} 60K15,45E99 }

\vskip 15mm
\section{\bf Introduction}

Many stochastic systems can be described by the abstract
mathematical model in the Banach space $\mathcal{B}(\mathbb{R}^d)$
of functions $\varphi(u), u \in \mathbb{R}^d$, called random
evolution model, introduced by Griego and Hersh {\rm \cite{GH, H1,
H2}}.

Asymptotic methods in the theory of random evolutions were applied
by many mathematicians (see, e.g. {\rm \cite{HP, KoPeTu, P}}).
Application of these methods to different stochastic systems may be
found in {\rm \cite{Kush}}. In {\rm \cite{Pin}} kinetic theory of
gases, isotropic transport on manifolds, stability of random
oscillators were studied by similar methods.

Among recent works in this area with applications in mathematical
biology we may remind results of Hillen and Othmer {\rm \cite{Hill,
OH}} (see also {\rm \cite{P1}}).

Namely, transport equations are used in mathematical biology to
model movement and growth of populations. Certain bacteria show the
following movement pattern: periods of strait runs alternate with
periods of random rotations which leads to reorientation of the
cells. We may model this movement by a velocity jump process, which
leads to a transport equation.

Thus, the linear transport equation $$\frac{\partial}{\partial
t}p(x,v,t)+v\nabla p(x,v,t)=-\lambda p(x,v,t)+\int_V\lambda
T(v,v')p(x,v',t)dv'$$ in which $p(x,v,t)$ represents the density of
particles at spatial position $x\in R^n$ moving with velocity $v\in
V\subset R^n$ at time $t \geq0$ arises. The turning rate $\lambda$
may be space- or velocity-dependent. The turning kernel or turn
angle distribution $T(v, v')$ gives the probability of a velocity
jump from $v'$ to $v$ if a jump occurs.

The evolutionary equation studied in our work generalizes the
transport equation described above. Analogical generalization of
telegraph-type equation is described in {\rm \cite{Sam}}.

One more example of application of asymptotic methods is described
in the works of Yin and Zhang {\rm \cite{YZ}}.

They study a model for production planning of a failure-prone
manufacturing system, which consists of machines whose production
capacity is modeled by a Markov or semi-Markov chain. In a
large-scale system various components may change at different rates.
Thus, the system may be decomposed and the states of the chain may
be aggregated. The introduction of a small parameter $\varepsilon>0$
makes the system belong to the category of two-time-scale systems.

The system is described by the equations of the type
$$\frac{dp^{\varepsilon}(t)}{dt}=\frac{1}{\varepsilon}p^{\varepsilon}(t)Q(t).$$
Here $p^{\varepsilon}(t)$ is the probability distribution of Markov
or semi-Markov chain, $Q(t)$ is the corresponding generator. As a
rule, the fast changing  process $p^{\varepsilon}(t)$ in the
physical or manufacturing systems is hard to analyze. But using
limit properties, obtained with the help of asymptotic expansions,
they replace the process by its "average" in the system under
consideration.

The problems of minimization of discounted cost function and optimal
control are studied in {\rm \cite{YZ}} with the use of asymptotic
approximations. Markov decision problems, stochastic control of
dynamical systems, numerical methods for control and optimization
are also analyzed.

We study the following generalization of the stochastic systems
described above.

A semi-Markov random evolution (SMRE) is created by a solution of
the evolutionary equation in Euclidean space $\mathbb{R}^d, d\geq1$
$$du^{\varepsilon}(t)/dt=v(u^{\varepsilon}(t);{\ae}(t/\varepsilon)), t\geq 0$$
with $u^{\varepsilon}(t)\in \mathbb{R}^d, v$ a given continuous
function from $\mathbb{R}^d\times E$ into $\mathbb{R}$,
$(E,\mathcal{E})$ a standard (polish) space. ${\ae}$ is a given
semi-Markov switching process ${\ae}(t), t\geq0$ on
$(E,\mathcal{E})$ given in terms of the semi-Markov kernel {\rm
\cite{Kor1}}
$$Q(x,B,t)=P(x,B)F_x(t), x\in E, B\in \mathcal{E}, t\geq
0$$ that defines the transition probabilities of a Markov renewal
process $(\ae_n, \tau_n, n\geq 0)$:
$$Q(x,B,t)=P\{\ae_{n+1}\in B, \theta_{n+1}\leq
t|\ae_n=x\}=$$ $$=P\{\ae_{n+1}\in B|\ae_n=x\}P\{\theta_{n+1}\leq
t|\ae_n=x\},$$ with $\theta_{n}:=\tau_{n+1}-\tau_n.$

The stochastic kernel
$$P(x,B)=P\{\ae_{n+1}\in B|\ae_n=x\}$$ defines the transition
probabilities of the embedded Markov chain $\ae_{n}={\ae}(\tau_n), $
$n\geq 0;$ the renewal moments are determined by the distribution
function of the sojourn times:
$$F_x(t)=P\{\theta_{n+1}\leq t|\ae_n=x\}=:P\{\theta_{x}\leq
t\}, x\in E.$$

We assume that the distribution functions $F_x(t), x\in E$ satisfy
the Cramer condition uniformly on $x\in E$: $$\sup\limits_{x\in
E}\int_0^{\infty} e^{ht}F_x(dt)\leq H<\infty, \eqno(1)$$ for all $
h>0.$

We denote by $Q$ the generator of the associated Markov process
${\ae}^0(t), t\geq0$ given by:
$$Q=q(x)(P-I), \eqno(2)$$ where the operator of transition probabilities $P$ is defined by
$$Pf(x)=\int_EP(x,dy)f(y), x\in E, \eqno(3)$$ for all bounded measurable real
valued $f$ defined on $E.$ $q(x)$ is defined by $$q(x):=1/m_1(x),$$
with $m_k(x)=\int_0^{\infty}s^k F_x(ds).$

We will see later that the equation for the regular part of a random
evolution is defined by the generator (2) of a uniformly ergodic
associated Markov process ${\ae}^0(t), t\geq0.$ The Banach space
$\mathcal{B}(E)$ is splitted as direct sum of the two subspaces {\rm
\cite{KoTu}}:
$$\mathcal{B}(E)=N_Q\bigoplus R_Q,$$ where
$N_Q:=\{\varphi:Q\varphi=0\}$ is the null-space of $Q$, and
$R_Q:=\{\psi:Q\varphi=\psi\}$ is the range of $Q$.

We define the projector $\Pi:$ $N_Q:=\Pi \mathcal{B}(E),
R_Q:=(I-\Pi)\mathcal{B}(E);
\Pi\varphi(x):=\widehat{\varphi}\mathbf{1},
\widehat{\varphi}:=\int_E \varphi(x)\pi(dx),$ where the stationary
distribution $\pi(B), B\in \mathcal{E}$ of the semi-Markov process
${\ae}(t),t\geq0$ satisfies the relations {\rm \cite{KoKo, KoTurb}}
$$\pi(dx)=\rho(dx)m_1(x)/\widehat{m},$$
$$\widehat{m}:=\int_E m_1(x)\rho(dx).$$ $\rho(B), B\in \mathcal{E}$
is by assumption the stationary distribution of the embedded Markov
chain ${\ae}_n, n\geq
 0,$ given by the formula $$\rho(B)=\int_E P(x,B)\rho(dx), \rho(E)=1.$$

Let us consider the Banach space $\mathcal{B}(\mathbb{R}^d)$ of
real-valued test-functions $\varphi(u), u\in \mathbb{R}^d$ which are
bounded with all their derivatives. We equip
$\mathcal{B}(\mathbb{R}^d)$ with $\sup$-norm
$$||\varphi||:=\sup\limits_{u\in
\mathbb{R}^d}|\varphi(u)|<C_{\varphi}$$ for some $C_{\varphi}>0$.

The random evolution in $\mathcal{B}(\mathbb{R}^d)$ is given by the
relation
$$\Phi_t^{\varepsilon}:=\varphi(u^{\varepsilon}(t)), \varphi\in \mathcal{B}(\mathbb{R}^d).\eqno(4)$$

In the present paper we shall investigate the asymptotic behavior of
SMRE (4) as $\varepsilon\to 0$ under the assumption of uniformly
ergodicity of the semi-Markov switching process ${\ae}(t)$ described
above and under the assumption of the existence of a global solution
of the deterministic equations
$$du_x(t)/dt=v(u_x(t);x), x\in E. \eqno(5)$$

Let us first consider the deterministic evolution
$$\Phi_x(t,u)=\varphi (u_x(t)), u_x(0)=u \in \mathbb{R}^d.$$ It generates a
corresponding semigroup on $\mathcal{B}(\mathbb{R}^d):$
$$\mathbb{V}_t(x)\varphi(u):=\varphi (u_x(t)), u_x(0)=u,$$ and its
generator has the form:
$$\mathbb{V}(x)\varphi(u)=v(u;x)\varphi '(u):=\sum_{k=1}^dv_k(u;x)\varphi_k '(u), \varphi(u)\in C^{\infty}(\mathbb{R}^d). \eqno(6)$$

Under the assumptions of [10] by the average principle the weak
convergence
$$u^{\varepsilon}(t)\Rightarrow \widehat{u} (t), \varepsilon\to 0 \eqno(7)$$
takes place. The average limit evolution $\widehat{u} (t), t\geq 0$
is defined by a solution of the average equation $$d\widehat{u}
(t)/dt=\widehat{v}(\widehat{u} (t)).$$

The average velocity $\widehat{v}(u), u\in \mathbb{R}^d$ is defined
by $$\widehat{v}(u)=\int_Ev(u;x)\pi(dx)$$ (i.e. by the average of
the initial velocity $v(u;x)$ over the stationary distribution
$\pi(B), B\in \mathcal{E}$).

The rate of convergence in (7) can be investigated in two
directions: \\ i) asymptotic analysis as $\varepsilon\downarrow 0$
of the fluctuations
$$\zeta^{\varepsilon}(t)=u^{\varepsilon}(t)-\widehat{u}(t); \eqno(8)$$
\\ ii) asymptotic analysis of the average deterministic evolution
$$\Phi_t^{\varepsilon}(u,x)=E[\varphi(u^{\varepsilon}(t))|u^{\varepsilon}(0)=u, {\ae}^{\varepsilon}(0)=x]. \eqno(9)$$

The asymptotic analysis of fluctuations (8) leads to the diffusion
approximation of the random evolution {\rm \cite{KoLim, Skor}}.

The asymptotic analysis of the evolution (9) is realized in what
follows by constructing the asymptotic expansion in terms of the
small parameter $\varepsilon$ $(\varepsilon\downarrow 0)$ in the
following form ($\tau=t/\varepsilon$):
$$\Phi_t^{\varepsilon}(u,x)=U_0(t)+\sum_{k=1}^{\infty}\varepsilon^k[U_k(t)+W_k(\tau)]. \eqno(10)$$

The asymptotic expansion (10) contains two parts: \\ i) the regular
term $U^{\varepsilon}(t):=U_0(t)+\sum_{k=1}^{\infty}\varepsilon^k
U_k(t),$ \\ ii)the singular term (boundary layer)
$W^{\varepsilon}(\tau):=\sum_{k=1}^{\infty}\varepsilon^k W_k(\tau),
\tau=t/\varepsilon.$

In addition the 'regular degeneration' (see {\rm \cite{Kor2, VaBu}}) of
the initial and 'boundary' condition
$$\begin{array}{c}
U_0(0)=\varphi(u)\mathbf{1}, U_k(0)+W_k(0)=0, k\geq 1, \\
W^{\varepsilon}(t)+U^{\varepsilon}(\varepsilon t)-\varphi(u)=0,
t\leq 0 \\
\end{array} \eqno(11)$$ has to be valid for all $x\in E, u\in \mathbb{R}^d$.

\begin{rem} The singular terms in (10) are determined by the boundary condition at infinity:
$$W_k(+\infty)=0, k\geq 1.$$
\end{rem}

Asymptotic expansions with 'boundary layers' were studied by many
authors (see {\rm \cite{Kor2, VaBu}}). In particular, functionals of
Markov and semi-Markov processes are investigated from this point of
view in {\rm \cite{KoPeTu, Sam, Tad}}.

In this work the asymptotic expansion (10) for the average
semi-Markov random evolution (9) is constructed by using the
solution of an integral Markov renewal equation. The algorithm for
the construction of explicit regular and singular terms and boundary
conditions is given in the following Theorem, which will be proved
in some stages (see the following Lemmas).

We use the following notations: let
$$L_kU(t):=\sum_{n=0}^k(-1)^n\mathbb{V}^n(x)PU^{(k-n)}(t), \eqno(12)$$ where
$U(t)$ is any smooth function, $U^{(k)}(t):=\frac{\partial^k
U(t)}{\partial t^k}$, $P$ is defined by (3), $\mathbb{V}(x)$ is
defined by (6),
$$\mu_k(x):=\frac{m_k(x)}{k!m_1(x)}, \mu_1(x):=1.$$

We also set:
$$\mathbf{Q}W(\tau):=\int_0^{\infty}F_x(ds)PW(\tau-s),$$
$$\psi^k(\tau):=\overline{F}^{(k)}(\tau)\mathbb{V}^kP\varphi(u),\psi^k_0(\tau):=\sum_{r=1}^{k-1}\mathbf{Q}^r
W_{k-r}(\tau),$$
$$\overline{F}^{(k)}(\tau):=\int_{\tau}^{\infty}\frac{s^{k-1}}{(k-1)!}\overline{F}_x(s)ds,
\mathbf{Q}^rW(\tau):=\int_0^{\infty}\frac{s^r}{r!}F_x(ds)\mathbb{V}^rPW(\tau-s),$$
$$\overline{F}_x(t):=1-F_x(t), \tau=t/\varepsilon.$$

\begin{thm} Under the conditions of uniform ergodicity of the underlying semi-Markov process and
the existence of a global solution in $\mathcal{B}(\mathbb{R}^d)$ of
the system (5), the asymptotic expansion of the semi-Markov
evolution
$$\Phi_t^{\varepsilon}(u,x)=E[\varphi(u^{\varepsilon}(t))|
u^{\varepsilon}(0)=u,{\ae}(0)=x]$$ has the form
$$\Phi_t^{\varepsilon}(u,x)=U^{\varepsilon}(t)+W^{\varepsilon}(\tau)=U_0(t)+\sum_{k=1}^{\infty}
\varepsilon^k(U_k(t)+W_k(\tau)),$$ where
$$U_0(t)=c_0(t)\mathbf{1}.$$ The function $ c_0(t)$ satisfies the
equation
$$\hat{v}(u)\frac{\partial c_0(t)}{\partial u}-\frac{\partial
c_0(t)}{\partial t}=0,$$ or
$$\partial c_0(t)/\partial t=\widehat{\mathbb{V}}c_0(t),$$ where (correspondingly to (6)):
$$\widehat{\mathbb{V}}\varphi(u):=\widehat{v}(u)\varphi'(u),$$
$\hat{v}(u):=\int_Ev(u;x)\pi(dx)$ is the average of the initial
velocity $v(u;x)$ over the stationary distribution of the
semi-Markov process. The initial condition is
$$c_0(0)=\varphi(u).$$

The regular terms are the following:
$$U_k(t)=\mathbb{R}_0\left(\sum_{n=1}^k\mu_n(x)L_nU_{k-n}(t)\right)+c_k(t),$$
where, according to {\rm \cite{KoTu}},
$\mathbb{R}_0=[Q+\Pi]^{-1}-\Pi.$

The functions $c_k(t)$ satisfy the equations
$$\hat{L}_1c_k(t)=-\Pi\mathfrak{L}_kc_0(t)-\ldots-\Pi\mathfrak{L}_1c_{k-1}(t), k\geq 1$$
where
$\Pi\mathfrak{L}_k:=\sum_{n=1}^k\Pi\mu_n(x)L_n\mathbb{R}_0\mathfrak{L}_{k-n}+\Pi\mu_{k+1}(x)L_{k+1},
\Pi\mathfrak{L}_0:=\Pi L_1=\hat{L}_1.$

The singular terms satisfy the Markov renewal equations $(k\geq 1)$:
$$\int_0^{\tau}F_x(ds)PW_k(\tau-s)-W_k(\tau)=\psi^k(\tau)-\psi^k_0(\tau)-\psi^k_1(\tau), \tau\geq
0$$ with
$$\psi^k_1(\tau):=\int_{\tau}^{\infty}F_x(ds)PW_k(\tau-s),$$
($P$ is given by (3)) and may be explicitly written in the form:
$$\begin{array}{c}
W_1(\tau)= \mathbf{R_0}[\psi^1(\tau)+\overline{F}_x(\tau)PU_1(0)+\int_{\tau}^{\infty}(\tau-s)F_x(ds)PU'_0(0)],\\
W_k(\tau)= \mathbf{R_0}[\psi^k(\tau)-\psi^k_0(\tau)+\overline{F}_x(\tau)PU_k(0)+\sum_{n=1}^k\int_{\tau}^{\infty}(\tau-s)^nF_x(ds)PU^{(n)}_{k-n}(0)],\\
  \end{array} $$ where
$\mathbf{R_0}$ is the Markov renewal operator {\rm \cite{Shu}}.

The initial conditions are given by: $$(I-\Pi)[U_k(0)+W_k(0)]=0,$$
$$ c_k(0)=-\Pi W_k(0),$$ $$U_k(0)=\left[\sum_{r=0}^{k-1}\int\pi(dx)\nu_{k-r}(x)L_{k-r}U_r(0)-\right.$$ $$\left.-\sum_{r=1}^{k-1}\int\rho(dx)\int_0^{\infty}\int_0^{\tau}\frac{s^r}{r!}F_x(ds)\mathbb{V}^r(x)PW_{k-r}(\tau-s)d\tau\right]/\widehat{m},$$
where $\nu_k(x)=(-1)^k[m_k(x)-\mu_{k+1}(x)].$
\end{thm}
{\it Proof.} The proof is a consequence of the Lemmas in sections
2-6.

\begin{rem} For sufficient conditions for the assumptions to hold see, e.g. [4,10].
The equation for $c_0(t)$ corresponds to the one given by the averaging theorem
{\rm \cite{KoLim}}. It states that the limit of the semi-Markov
evolution
$$\Phi^{\varepsilon}_t(u,x)\Rightarrow \widehat{\Phi}_t(u),
\varepsilon\to 0$$ satisfies the equation
$$\frac{\partial \widehat{\Phi}_t(u)}{\partial t}=\widehat{\mathbb{V}}\widehat{\Phi}_t(u).$$

Here $\widehat{\Phi}_t(u)$ is the deterministic evolution
$$\widehat{\Phi}_t(u)=\varphi (u(t)), u(0)=u.$$ The average limit evolution ${u} (t), t\geq 0$ is
defined by a solution of the average equation $$d{u}
(t)/dt=\widehat{v}({u} (t)).$$
\end{rem}

\section{\bf Markov renewal equation}

\begin{lem} The semi-Markov evolution $\Phi_t^{\varepsilon}(u,x)$
satisfies the equation
$$\int_0^{\infty}F_x(ds)\mathbb{V}_{\varepsilon
s}(x)P\Phi_{t-\varepsilon s}^{\varepsilon}(u,x)-
\Phi_{t}^{\varepsilon}(u,x)=$$
$$=\varepsilon\mathbb{V}(x)\int_{\tau}^{\infty}\overline{F}_x(s)
\mathbb{V}_{\varepsilon s}(x)\varphi(u)ds, \eqno(13)$$ where
$\tau=t/\varepsilon.$ \end{lem}

{\it Proof.} Using the first jump moment of the switching process
$\theta_x$, we have:
$$\Phi_t^{\varepsilon}(u,x)=E_{u,x}[\varphi(u^{\varepsilon}(t));\theta_x>
t/\varepsilon]+E_{u,x}[\varphi(u^{\varepsilon}(t));\theta_x\leq
t/\varepsilon]=$$ $$=
\overline{F}_x(t/\varepsilon)\mathbb{V}_t(x)P\varphi(u)+\int_0^{t/\varepsilon}
F_x(ds)\mathbb{V}_{\varepsilon s}(x)P\Phi_{t-\varepsilon
s}^{\varepsilon}(u,x).$$

So:
$$\Phi_{t}^{\varepsilon}(u,x)-\int_0^{t/\varepsilon}F_x(ds)\mathbb{V}_{\varepsilon
s}(x)P\Phi_{t-\varepsilon
s}^{\varepsilon}(u,x)=\overline{F}_x(\tau)\mathbb{V}_t(x)P\varphi(u).$$

Extending by the continuity, $\Phi_{t-\varepsilon
s}^{\varepsilon}(u,x)=\varphi(u), t-\varepsilon s\leq 0, $ let us rewrite
the latter equation in the form:
$$\Phi_{t}^{\varepsilon}(u,x)-\int_0^{\infty}F_x(ds)\mathbb{V}_{\varepsilon
s}(x)P\Phi_{t-\varepsilon s}^{\varepsilon}(u,x)=
\overline{F}_x(\tau)\mathbb{V}_t(x)P\varphi(u)-$$ $$-
\int_{\tau}^{\infty} F_x(ds)\mathbb{V}_{\varepsilon
s}(x)P\Phi_{t-\varepsilon
s}^{\varepsilon}(u,x)=\overline{F}_x(\tau)\mathbb{V}_t(x)P\varphi(u)-\int_{\tau}^{\infty}
F_x(ds)\mathbb{V}_{\varepsilon s}(x)P\varphi(u).$$

So, we have:
$$\Phi_{t}^{\varepsilon}(u,x)-\int_0^{\infty}F_x(ds)\mathbb{V}_{\varepsilon
s}(x)P\Phi_{t-\varepsilon
s}^{\varepsilon}(u,x)=\overline{F}_x(\tau)\mathbb{V}_t(x)P\varphi(u)-$$
$$- \overline{F}_x(s)\mathbb{V}_{\varepsilon
s}(x)P\varphi(u)|_{\tau}^{\infty}-\varepsilon\mathbb{V}(x)\int_{\tau}^{\infty}\overline{F}_x(s)
\mathbb{V}_{\varepsilon s}(x)\varphi(u)ds.$$

This gives equation (13). $\Box$

\section{\bf Equations for the regular terms}

\begin{lem} The equation for the regular part of the asymptotics has the form:
$$QU^{\varepsilon}(t)=\left[\sum_{k=1}^{\infty}\varepsilon^k\mu_k(x)L_k\right]U^{\varepsilon}(t).\eqno(14)$$
\end{lem}

{\it Proof.} We use the equality:
$$aPb-1=(P-1)+(a-1)P+P(b-1)+(a-1)P(b-1),$$ where
$a=\mathbb{V}_{\varepsilon
s}(x)=I+\sum_{k=1}^{\infty}\varepsilon^k\frac{s^k}{k!}\mathbb{V}^k(x),
b=\Phi_{t-\varepsilon s}^{\varepsilon}=\sum_{k=0}^{\infty}(-1)^k
\varepsilon^k\frac{s^k}{k!}\Phi^{(k)}_t(u,x).$

Let us rewrite (13) in the following way:
$$(P-I)\Phi^{\varepsilon}_t(u,x)+\int_0^{\infty}F_x(ds)\left(\sum_{k=1}^{\infty}\varepsilon^k\frac{s^k}{k!}\mathbb{V}^k(x)\right)
\times$$ $$\times
P\Phi^{\varepsilon}_t(u,x)+\int_0^{\infty}F_x(ds)P\left(\sum_{k=1}^{\infty}(-1)^k
\varepsilon^k\frac{s^k}{k!}\Phi^{(k)}_t(u,x)\right)+\int_0^{\infty}F_x(ds)
\left(\sum_{k=1}^{\infty}\varepsilon^k\frac{s^k}{k!}\times\right.$$
$$\left.\times\mathbb{V}^k(x)\right)P\left(\sum_{k=1}^{\infty}(-1)^k
\varepsilon^k\frac{s^k}{k!}\Phi^{(k)}_t(u,x)\right)=\varepsilon
\mathbb{V}(x)\int_{\tau}^{\infty}\overline{F}_x(s)\mathbb{V}_{\varepsilon
s}(x)P\varphi(u)ds.$$

Substituting (10) for the regular part, we have:
$$(P-I)U^{\varepsilon}(t)=-\int_0^{\infty}F_x(ds)\left(\sum_{k=1}^{\infty}\varepsilon^k\frac{s^k}{k!}\mathbb{V}^k(x)\right)
PU^{\varepsilon}(t)-\int_0^{\infty}F_x(ds)P\left(\sum_{k=1}^{\infty}(-1)^k\times\right.$$
$$\left.\times\varepsilon^k\frac{s^k}{k!}U^{\varepsilon (k)}(t)\right)-\int_0^{\infty}F_x(ds)
\left(\sum_{k=1}^{\infty}\varepsilon^k\frac{s^k}{k!}\mathbb{V}^k(x)\right)P\left(\sum_{k=1}^{\infty}(-1)^k
\varepsilon^k\frac{s^k}{k!}U^{\varepsilon (k)}(t)\right).$$

Gathering the terms with the same degree of $\varepsilon$, we obtain:
$$(P-I)U^{\varepsilon}(t)=\sum_{k=1}^{\infty}\varepsilon^k\left[-\int_0^{\infty}F_x(ds)\frac{s^k}{k!}\mathbb{V}^k(x)
PU^{\varepsilon}(t)-\int_0^{\infty}F_x(ds)\left(\sum_{n=1}^{k-1}(-1)^{n}
\times\right.\right.$$ $$\left.\left.\times
\frac{s^k}{(k-n)!n!}\mathbb{V}^{(k-n)}(x) PU^{\varepsilon
(n)}(t)\right)-\int_0^{\infty}(-1)^kF_x(ds)P\frac{s^k}{k!}U^{\varepsilon
(k)}(t)\right]=
\sum_{k=1}^{\infty}\varepsilon^k\frac{m_k(x)}{k!}L_kU^{\varepsilon}(t).
$$

To obtain (14) we should divide the last equality by $m_1(x)$.

Then the lemma is proved. $\Box$

If we put into (14) the expansion
$U^{\varepsilon}(t)=\sum_{k=0}^{\infty}\varepsilon^kU_k(t)$ and
gather together the terms with the same degree of $\varepsilon$, we
obtain the following corollary.

\begin{cor} The regular terms of the asymptotics satisfy the following system of equations:
$$\left\{\begin{array}{c}
QU_0(t)=0, \\
QU_k(t)=\sum_{n=1}^k\mu_n(x)L_nU_{k-n}(t), k\geq 1\\
\end{array}\right. \eqno(15)$$
\end{cor}

The first regular term $U_0(t)$ belongs to $N_Q,$ according to the
first equation of the system (15).

Hence:
$$U_0(t)=c_0(t)\mathbf{1},$$ where $ c_0(t)$ is a scalar function that does not depend on
$x$. To obtain the equation for $c_0(t)$ the solvability condition
for the following equation is used: $$\Pi L_1c_0(t)=0,$$ which leads
to the average equation $$\hat{L}_1c_0(t)=0, \hat{L}_1 \Pi=\Pi
L_1\Pi. \eqno (16)$$

\begin{cor} The function $c_0(t)$ satisfies the equation with initial condition:
$$\hat{v}(u)\frac{\partial c_0(t)}{\partial u}-\frac{\partial
c_0(t)}{\partial t}=0,$$
$$c_0(0)=\varphi(u).$$
\end{cor}

{\it Proof.} The explicit form of $\hat{L}_1$ may easily be found
using (12):
$$\hat{L}_1c_0(t):=\hat{v}(u)\frac{\partial c_0(t)}{\partial
u}-\frac{\partial c_0(t)}{\partial t}, \hat{v}(u):=\Pi
v(u,x):=\int_Ev(u,x)\pi(dx).$$

From this and (16) we obtain the equation in corollary 3.3. The
initial condition follows from the formulae (11).

The corollary is thus proved. $\Box$

For $U_1(t)$ we obtain:
$$U_1(t)=\mathbb{R}_0L_1U_0(t)+c_1(t)=\mathbb{R}_0L_1c_0(t)+c_1(t).$$

Using the solvability condition for the third equation of the system
(15), we have:
$$\Pi L_1\mathbb{R}_0L_1c_0(t)+\Pi
\mu_2(x)L_2c_0(t)+\hat{L}_1c_1(t)=0,$$ or
$$\hat{L}_1c_1(t)=-\Pi\mathfrak{L}_1c_0(t),$$ where $\Pi\mathfrak{L}_1:=\Pi L_1\mathbb{R}_0L_1+\Pi\mu_2(x) L_2.$

For $U_k(t)$ we obtain by analogy:
$$U_k(t)=\mathbb{R}_0\left(\sum_{n=1}^k\mu_n(x)L_nU_{k-n}(t)\right)+c_k(t),$$

$$\hat{L}_1c_k(t)=-\Pi\mathfrak{L}_kc_0(t)-\ldots-\Pi\mathfrak{L}_1c_{k-1}(t),$$
where
$\Pi\mathfrak{L}_k:=\sum_{n=1}^k\Pi\mu_n(x)L_n\mathbb{R}_0\mathfrak{L}_{k-n}+\Pi\mu_{k+1}(x)L_{k+1},
\Pi\mathfrak{L}_0:=\Pi L_1.$

The regular terms $U_k(t), k\geq 1$ contain two parts: one in the
null-space $c_k(t)$ and one in the space of values $R_Q$:
$$U_k(t)=c_k(t)\mathbf{1}+U_k^R(t), U_k^R(t)\in R_Q.$$

The scalar functions $c_k(t), k\geq1$ satisfy the same equations
with the operator $\widehat{L}_1$ but with different right-hand
side. Hence the unique solution can be defined by using the initial
condition $c_k(0).$ Meanwhile the second part $U_K^R(t)$ is defined
uniquely by the relation
$$U_K^R(t)=\mathbb{R}_0\left(\sum_{n=1}^k\mu_n(x)L_nU_{k-n}(t)\right).$$

It is worth noticing that the regular part $U_K^R(t)$ in $R_Q$ is
determined uniquely for $t$ on the entire real line
$\mathbb{R}=(-\infty, +\infty)$.

\section{\bf Equations for the singular terms}

The equation for the singular term $W(\tau), \tau=t/\varepsilon
(\varepsilon>0)$ is constructed in a very simple way, because the
shift by $-\varepsilon s$ may be transformed into the convolution
$$\mathcal{L}^{\varepsilon}W^{\varepsilon}(\tau)=\int_0^{\infty}F_x(ds)\mathbb{V}_{\varepsilon s}(x)PW^{\varepsilon}(\tau-s)-W^{\varepsilon}(\tau)=\varepsilon\psi_{\varepsilon}(\tau), \eqno(17)$$
where
$\varepsilon\psi_{\varepsilon}(\tau)=\varepsilon\mathbb{V}(x)\int_{\tau}^{\infty}\overline{F}_x(s)
\mathbb{V}_{\varepsilon s}(x)\varphi(u)ds.$

The algebraic identity $$VP-I=P-I+(V-I)P$$ provides the
representation
$$(\mathbf{Q}-I)W^{\varepsilon}(\tau)+\mathbf{Q}_1^{\varepsilon}W^{\varepsilon}(\tau)=\varepsilon\psi_{\varepsilon}(\tau),$$
where
$$(\mathbf{Q}-I)W^{\varepsilon}(\tau):=\int_0^{\infty}F_x(ds)PW^{\varepsilon}(\tau-s)-W^{\varepsilon}(\tau)$$
is the renewal operator on the real line $-\infty<\tau<+\infty$ and
$$\mathbf{Q}_1^{\varepsilon}W^{\varepsilon}(\tau):=\int_0^{\infty}F_x(ds)[\mathbb{V}_{\varepsilon s}(x)-I]PW^{\varepsilon}(\tau-s).$$

Now the equations for the singular terms are constructed in the
usual way:

\begin{lem} The equations for the singular terms have the following form:
$$\begin{array}{c}
(\mathbf{Q}-I)W_1(\tau)=\psi^1(\tau),\\
(\mathbf{Q}-I)W_k(\tau)=\psi^k(\tau)-\psi^k_0(\tau),\\
  \end{array} \eqno(18)$$
\end{lem}

{\it Proof.} If we put the following expansion
$W^{\varepsilon}(\tau)=\sum_{k=1}^{\infty}\varepsilon^kW_k(\tau)$
into equation (17) we obtain:
$$\int_0^{\infty}F_x(ds)\left[I+\sum_{k=1}^{\infty}\varepsilon^k\frac{s^k}{k!}\mathbb{V}^k(x)\right]
P\left[\sum_{k=1}^{\infty}\varepsilon^kW_k(\tau-s)\right]-$$
$$-\sum_{k=1}^{\infty}\varepsilon^kW_k(\tau)=
\int_{\tau}^{\infty}\overline{F}_x(s)\left[\sum_{k=1}^{\infty}\varepsilon^k\frac{s^{k-1}}{(k-1)!}
\mathbb{V}^k(x)\right]P\varphi(u)ds.$$

So, we have: $$\varepsilon
[\mathbf{Q}-I]W_1(\tau)+\sum_{k=2}^{\infty}\varepsilon^k[\mathbf{Q}-I]W_k(\tau)+\sum_{k=2}^{\infty}\varepsilon^k
\sum_{r=1}^{k-1}\mathbf{Q}^r\times$$ $$\times
W_{k-r}(\tau)=\sum_{k=1}^{\infty}\varepsilon^k\overline{F}^{k}(\tau)
\mathbb{V}^k(x)P\varphi(u),$$ and gathering the terms of the same
degree of $\varepsilon$, we obtain (18).

The lemma is thus proved. $\Box$

\begin{cor} The singular terms of the asymptotics have the following form:
$$\begin{array}{c}
W_1(\tau)= \mathbf{R_0}[\psi^1(\tau)-\psi^1_1(\tau)],\\
W_k(\tau)= \mathbf{R_0}[\psi^k(\tau)-\psi^k_0(\tau)-\psi^k_1(\tau)], k\geq 2.\\
  \end{array} $$ Here $\mathbf{R}_0$ is the Markov renewal operator {\rm \cite{Shu}}.
\end{cor}

{\it Proof.} The extended Markov renewal equation (18) for the
singular terms may be transformed into a standard form as follows
$$\int_0^{\tau}F_x(ds)PW_k(\tau-s)-W_k(\tau)=\psi^k(\tau)-\psi^k_0(\tau)-\psi^k_1(\tau), \tau\geq
0$$ with
$$\psi^k_1(\tau):=\int_{\tau}^{\infty}F_x(ds)PW_k(\tau-s).$$

Following [9] we may find the solution of standard Markov renewal
equation using the Markov renewal operator $\mathbf{R_0}$.

Hence the corollary is proved. $\Box$

From the corollary we have that the singular term can be defined
after continuation to the negative real line for $\tau<0.$ The
additional boundary condition
$$W_k(+\infty)=0$$ and the 'regular degeneration' of boundary conditions (11)
 provide the unique determination of the singular term.

\section{\bf Regularity of boundary conditions}

The initial extended Markov renewal equation (13) takes place under
the continuously differentiable extension of the regular part
$U^{\varepsilon}(t)$ on the negative real line for $t<0.$ That is
the concordance of the boundary conditions (11) have to be valid, or
in other form we should have
$$W^{\varepsilon}(t)+U^{\varepsilon}(\varepsilon t)+U_0(\varepsilon
t)=0, t<0. \eqno(19)$$

The representations (19) and (11) and the Taylor formula for
$\tau<0$:
$$\varphi(u)=\Phi_{\varepsilon\tau}(u,x)|_{\tau<0}=U_0(0)+\sum_{k=1}^{\infty}\varepsilon^k
\frac{\tau^k}{k!}U^{(k)}_0(0)+\varepsilon U_1(0)+$$
$$+\varepsilon\sum_{k=1}^{\infty}\varepsilon^k
\frac{\tau^k}{k!}U^{(k)}_1(0)+\ldots+\sum_{k=1}^{\infty}\varepsilon^kW_k(\tau)$$
give the boundary conditions for the singular terms in the following
form
$$W_k(\tau)=W_k(0)-\sum_{n=1}^{k}\frac{\tau^n}{n!}U^{(n)}_{k-n}(0), \tau<0, k\geq1
\eqno(20)$$ and the additional initial condition $$W_k(0)+U_k(0)=0,
k\geq1. \eqno(21)$$

The condition (20) may be used to construct a solution of the Markov
renewal equation (18) for the singular terms. The boundary condition
for singular terms $W_k(+\infty)=0, k\geq1$ and the Markov renewal
limit theorem {\rm \cite{Shu}} provide by a calculation: $$\Pi
W_k(0)\in N_Q$$ and, as a result, we obtain an initial condition for
the regular part in $N_Q$
$$c_k(0)=-\Pi W_k(0).$$

The last step to verify the regularity of boundary conditions (11)
is to establish the following relations: $$[P-I](W_k(0)+U_k(0))=0,
k\geq1,$$ which are, in fact, equivalent to (11).

These are the regular and corresponding singular terms in the
subspace of values $R_Q$ which compensate each other without any
further assumption. A good formulation of the problem of asymptotic
expansion leads to the regularity of the boundary conditions (11).

\begin{lem} The regular degeneration equation has the form:
$$(P-I)[U^{\varepsilon}(0)+W^{\varepsilon}(0)]=0.$$
\end{lem}

{\it Proof.} Let us consider the equation (18). For $W_1(\tau)$ we
have:
$$\int_0^{\infty}F_x(ds)PW_1(\tau-s)-W_1(\tau)=\int_{\tau}^{\infty}\overline{F}_x(ds)\mathbb{V}(x)PU_0(0).$$

For $\tau=0$ we obtain:
$$\left[\int_0^{\infty}F_x(ds)PW_1(0)-W_1(0)\right]+\int_0^{\infty}F_x(ds)P(W_1(-s)-W_1(0))=\int_{0}^{\infty}\overline{F}_x(ds)\mathbb{V}(x)PU_0(0).$$

Using the equality (20) we have:
$$(P-I)W_1(0)=-\int_0^{\infty}sF_x(ds)PU'_0(0)+\int_{0}^{\infty}s{F}_x(ds)\mathbb{V}(x)PU_0(0). \eqno(22)$$

For the corresponding regular term, we get from (14) and (12):
$$(P-I)U_1(0)=m_1(x)PU'_0(0)-m_1(x)\mathbb{V}(x)PU_0(0).$$

So, the following equality is true:
$$(P-I)\left[U_1(0)+W_1(0)\right]=0.$$

By induction, if we have:
$$(P-I)\left[U_n(0)+W_n(0)\right]=0, n=1,\ldots,k, \eqno(23)$$ then for
$U_{k+1}(0)$ and $W_{k+1}(0)$ we obtain:
$$\int_0^{\infty}F_x(ds)PW_{k+1}(\tau-s)-W_{k+1}(\tau)=\int_{\tau}^{\infty}\overline{F}_x(ds)\sum_{n=1}^{k+1}\frac{s^{n-1}}{(n-1)!}\mathbb{V}^n(x)PU_0(0)-$$
$$-\int_0^{\infty}F_x(ds)\sum_{n=1}^k\frac{s^n}{n!}\mathbb{V}^n(x)PW_{k-n+1}(\tau-s).$$

For $\tau=0$ we have:
$$\left[\int_0^{\infty}F_x(ds)PW_{k+1}(0)-W_{k+1}(0)\right]+\int_0^{\infty}F_x(ds)P(W_{k+1}(-s)-W_{k+1}(0))=$$ $$=\int_{0}^{\infty}\overline{F}_x(ds)\frac{s^{k}}{k!}\mathbb{V}^{k+1}(x)PU_0(0)
-\int_0^{\infty}F_x(ds)\sum_{n=1}^k\frac{s^n}{n!}\mathbb{V}^n(x)PW_{k-n+1}(-s).$$

By the equality (20) we obtain:
$$(P-I)W_{k+1}(0)=\int_0^{\infty}F_x(ds)P\sum_{n=1}^{k+1}\frac{(-s)^n}{n!}U_{k-n+1}^{(n)}(0)+\int_{0}^{\infty}{F}_x(ds)\frac{s^{k+1}}{(k+1)!}\mathbb{V}^{k+1}(x)PU_0(0)-$$
$$-\int_0^{\infty}F_x(ds)\sum_{n=1}^k\frac{s^n}{n!}\mathbb{V}^n(x)P\left[W_{k-n+1}(0)-\sum_{n=1}^{k-n+1}\frac{(-s)^n}{n!}U_{k-n+1}^{(n)}(0)\right].$$

By induction, as soon as (21) and (23) are true, we may write:
$$(P-I)W_{k+1}(0)=-\int_0^{\infty}F_x(ds)P\sum_{n=1}^{k+1}\frac{(-s)^n}{n!}U_{k-n+1}^{(n)}(0)+\int_{0}^{\infty}{F}_x(ds)\frac{s^{k+1}}{(k+1)!}\mathbb{V}^{k+1}(x)PU_0(0)-$$
$$-\int_0^{\infty}F_x(ds)\sum_{n=1}^k\frac{s^n}{n!}\mathbb{V}^n(x)P\left[-U_{k-n+1}(0)-\sum_{n=1}^{k-n+1}\frac{(-s)^n}{n!}U_{k-n+1}^{(n)}(0)\right]. \eqno(24)$$

For the corresponding regular term, we have from (14) and (12):
$$(P-I)U_{k+1}(0)=-\int_0^{\infty}F_x(ds)P\sum_{n=1}^{k+1}\frac{(-s)^n}{n!}U^{(n)}_{k-n+1}(0)-\int_{0}^{\infty}{F}_x(ds)\sum_{n=1}^{k+1}\frac{s^n}{n!}\mathbb{V}^n(x)PU_{k-n+1}(0)-$$
$$-\int_0^{\infty}F_x(ds)\sum_{n=1}^k\frac{s^n}{n!}\mathbb{V}^n(x)P\sum_{r=1}^{k-n+1}\frac{(-s)^r}{r!}U_{k-n-r+1}^{(r)}(0).$$

It is easy to see that the following equality is true:
$$(P-I)\left[U_{k+1}(0)+W_{k+1}(0)\right]=0.$$

So, by induction the lemma is proved. $\Box$

\begin{cor}
$$(I-P)[U^{\varepsilon}(0)+W^{\varepsilon}(0)]=0,$$ or, in other words,
$$(I-\Pi)[U_k(0)+W_k(0)]=0.$$
\end{cor}

The proof is obvious.

So, we can see that in the space of values of the operator $Q$ the
regular and singular parts of the solution fulfil the initial
conditions (11).

At the same time, in the null-space of $Q$ the initial conditions for the
regular terms are determined by the initial conditions for the singular
terms, so we have

\begin{cor} $$c_k(0)=-\Pi W_k(0), k\geq 1.$$
\end{cor}
{\it Proof.} We obtain obviously: $\Pi[W_k(0)+U_k(0)]=\Pi
W_k(0)+c_k(0)=0.$ $\Box$

\begin{cor} The singular terms of the asymptotics have the following explicit
form:
$$\begin{array}{c}
W_1(\tau)= \mathbf{R_0}[\psi^1(\tau)+\overline{F}_x(\tau)PU_1(0)+\int_{\tau}^{\infty}(\tau-s)F_x(ds)PU'_0(0)],\\
W_k(\tau)= \mathbf{R_0}[\psi^k(\tau)-\psi^k_0(\tau)+\overline{F}_x(\tau)PU_k(0)+\sum_{n=1}^k\int_{\tau}^{\infty}(\tau-s)^nF_x(ds)PU^{(n)}_{k-n}(0)],\\
  \end{array} $$
\end{cor}

{\it Proof.} Using formulae (20) and corollary 5.3 we easily obtain:
$$\int_{\tau}^{\infty}F_x(ds)PW_k(\tau-s)=\int_{\tau}^{\infty}F_x(ds)P[-U_k(0)-
\sum_{n=1}^k(\tau-s)^nU^{(n)}_{k-n}(0)]=$$
$$=-\overline{F}_x(\tau)PU_k(0)-
\sum_{n=1}^k\int_{\tau}^{\infty}(\tau-s)^nF_x(ds)PU^{(n)}_{k-n}(0).$$

This proves the corollary. $\Box$

\section{\bf Initial conditions for the regular terms}

We are going to write down an algorithm for the construction of
initial conditions (at $t=0$) for the regular terms using the
boundary conditions for the singular terms as $\tau\to\infty$ (see
Remark 1.1). For the first singular term $W_1(\tau)$ we have the
equation (see (18)):
$$\int_0^{\infty}Q(ds)W_1(\tau-s)-W_1(\tau)=\overline{F}^{(1)}_x(\tau)\mathbb{V}(x)P\varphi(u), \eqno(25)$$
where
$\overline{F}^{(1)}_x(\tau)=\int_{\tau}^{\infty}\overline{F}_x(s)ds.$

Separating the first integral into two parts, we obtain:
$$\int_0^{\tau}Q(ds)W_1(\tau-s)-W_1(\tau)=\overline{F}^{(1)}(\tau)\mathbb{V}(x)P\varphi(u)-\int_{\tau}^{\infty}Q(ds)W_1(\tau-s).$$

According to the renewal theorem {\rm \cite{Shu}} we have for
$\tau\to\infty$:
$$0=W_1(\infty)=\left(\int\rho(dx)\int_0^{\infty}\int_{\tau}^{\infty}\overline{F}_x(s)dsd\tau \mathbb{V}(x)P\varphi(u)
-\right.$$
$$\left.-\int\rho(dx)\int_0^{\infty}\int_{\tau}^{\infty}Q(ds)W_1(\tau-s)d\tau\right)/\widehat{m}, \eqno(26)$$
where $\hat{m}=\int\rho(dx)m_1(x).$

For $\tau<0$ we have from (20): $$W_1(\tau)=W_1(0)-\tau U'_0(0).
\eqno(27)$$



Substituting the correlation (27) into equation (26), we obtain:
$$0=\left(\int\rho(dx)\int_0^{\infty}\int_{\tau}^{\infty}\overline{F}_x(s)dsd\tau
\mathbb{V}(x)P\varphi(u) -\right.$$
$$\left.-\int\rho(dx)\int_0^{\infty}\int_{\tau}^{\infty}Q(ds)[W_1(0)-(\tau-s)U'_0(0)]d\tau\right)/\widehat{m}=$$
$$=\left(\int\rho(dx)\int_0^{\infty}\overline{F}_x^{(1)}(s)d\tau
\mathbb{V}(x)P\varphi(u) -\right.$$
$$\left.-\int\rho(dx)\int_0^{\infty}\int_{\tau}^{\infty}F_x(s)PW_1(0)dsd\tau
+\int\rho(dx)\int_0^{\infty}\int_{\tau}^{\infty}Q(ds)(\tau-s)U'_0(0)d\tau\right)/\widehat{m}=$$
$$=\left(\int\rho(dx)\frac{m_2(x)}{2}
\mathbb{V}(x)P\varphi(u) -\right.$$
$$\left.-\int\rho(dx)m_1(x)PW_1(0)
-\int\rho(dx)\frac{m_2(x)}{2}PU'_0(0)\right)/\widehat{m}=$$
$$=\left(-\int\rho(dx)m_1(x)\mu_2(x)L_1(x)\varphi(u)-\int\rho(dx)m_1(x)(P-I)W_1(0)-\right.$$
$$\left.-\int\rho(dx)m_1(x)W_1(0)\right)/\widehat{m}=$$ $$=\left(-\int\rho(dx)m_1(x)\mu_2(x)L_1(x)\varphi(u)-\int\rho(dx)m_1(x)(P-I)W_1(0)\right)/\widehat{m}-\Pi W_1(0)=
$$
$$=\left(-\int\rho(dx)m_1(x)\mu_2(x)L_1(x)\varphi(u)-\int\rho(dx)m_1(x)(P-I)W_1(0)\right)/\widehat{m}+c_1(0),
\eqno(28)$$ here $\mu_2(x)=\frac{m_2(x)}{2m_1(x)}.$

We may rewrite (22) in the form:
$$[P-I]W_1(0)=m_1(x)[\mathbb{V}(x)P\varphi(u)-PU'_0(0)]=L_1(x)\varphi (u).$$

If we put the last equality into (28), we have finally:
$$0=\left(-\int\rho(dx)m_1(x)\mu_2(x)L_1(x)\varphi(u)+\int\rho(dx)m_1^2(x)L_1(x)\varphi
(u)\right)/\widehat{m}+c_1(0),$$ or

$$c_1(0)=\Pi U_1(0)=\int\pi(dx)\nu_1(x)L_1(x)\varphi(u)/\widehat{m},$$
where $\pi(dx)=\rho(dx)m_1(x),
\nu_1(x)=\mu_2(x)-m_1(x)=\frac{m_2(x)-2m_1^2(x)}{2m_1(x)}.$

\begin{rem} It is easily seen that $\nu_1(x)=0$, when $F_x(t)$ is
distributed exponentially. In this case, we obviously have:
$$c_1(0)=\Pi U_1(0)=0.$$

Note, that $\nu_1(x)$ may be either positive or negative, see, e.g
{\rm \cite{KoLim}}.
\end{rem}

We shall describe the algorithm for computing the next terms of the
asymptotics in the case of $W_2(\tau)$:
$$\int_0^{\infty}Q(ds)W_2(\tau-s)-W_2(\tau)=\overline{F}^{(2)}_x(\tau)\mathbb{V}^2(x)P\varphi(u)-
\int_0^{\infty}\frac{s}{1!}F_x(ds)\mathbb{V}(x)PW_1(\tau-s),
\eqno(29)$$ where
$\overline{F}^{(2)}(\tau)=\int_{\tau}^{\infty}s\overline{F}_x(s)ds.$

Separating the first integral, we obtain:
$$\int_0^{\tau}Q(ds)W_2(\tau-s)-W_2(\tau)=\overline{F}^{(2)}(\tau)\mathbb{V}^2(x)P\varphi(u)-
\int_0^{\infty}\frac{s}{1!}F_x(ds)\mathbb{V}(x)PW_1(\tau-s)-$$
$$-\int_{\tau}^{\infty}Q(ds)W_2(\tau-s).$$

According to the renewal theorem {\rm \cite{Shu}} we have for
$\tau\to\infty$:
$$0=W_2(\infty)=\left(\int\rho(dx)\int_0^{\infty}\int_{\tau}^{\infty}s\overline{F}_x(s)dsd\tau \mathbb{V}^2(x)P\varphi(u)
-\right.$$
$$-\int\rho(dx)\int_0^{\infty}\left[\int_0^{\tau}\frac{s}{1!}F_x(ds)\mathbb{V}(x)PW_1(\tau-s)d\tau+\right.$$ $$\left.\left.+\int_{\tau}^{\infty}\frac{s}{1!}F_x(ds)\mathbb{V}(x)PW_1(\tau-s)d\tau\right]-\int\rho(dx)\int_0^{\infty}\int_{\tau}^{\infty}Q(ds)W_2(\tau-s)d\tau\right)/\widehat{m}.
\eqno(30)$$

For $\tau<0$ we have from (20): $$W_2(\tau)=W_2(0)-\tau
U'_1(0)-\tau^2U_0''(0). \eqno(31)$$


Substituting (31) into (30), we obtain, like in the case of
$W_1(\tau)$:
$$0=\left(\int\rho(dx)\int_0^{\infty}\int_{\tau}^{\infty}s\overline{F}_x(s)dsd\tau
\mathbb{V}^2(x)P\varphi(u)-\right.$$ $$-
\int\rho(dx)\int_0^{\infty}\left[\int_0^{\tau}\frac{s}{1!}F_x(ds)\mathbb{V}(x)PW_1(\tau-s)d\tau+\int_{\tau}^{\infty}\frac{s}{1!}F_x(ds)\mathbb{V}(x)P\left\{W_1(0)-\right.\right.$$
$$-\left.\left.(\tau-s)U_0'(0)\right\}d\tau\right] -$$
$$\left.-\int\rho(dx)\int_0^{\infty}\int_{\tau}^{\infty}Q(ds)[W_2(0)-(\tau-s)
U'_1(0)-(\tau-s)^2U_0''(0)]d\tau\right)/\widehat{m}=$$
$$=\left(-\int\rho(dx)m_1(x)\mu_2(x)L_1(x)U_1(0)+\int\rho(dx)m_1(x)\mu_3(x)L_2(x)U_0(0)- \right.$$
$$-\int\rho(dx)m_1(x)(P-I)W_2(0)-$$ $$\left.-
\int\rho(dx)\int_0^{\infty}\int_0^{\tau}\frac{s}{1!}F_x(ds)\mathbb{V}(x)PW_1(\tau-s)d\tau
\right)/\widehat{m}+c_2(0). \eqno(32)$$

We obtain from (24) in case of $k=1$:
$$[P-I]W_2(0)=m_2(x)\mathbb{V}^2(x)P\varphi
(u)-m_1(x)PU'_1(0)+m_2(x)PU''_0(0)+$$
$$+m_1(x)\mathbb{V}(x)PU_1(0)-m_2(x)\mathbb{V}(x)PU'_0(0)
=m_2(x)L_2(x)U_0(0)-m_1(x)L_1(x)U_1(0).$$

Substituting the last equality into (32), we finally get:
$$c_2(0)=\left[\int\pi(dx)\nu_2(x)L_2(x)U_0(0)+
\int\pi(dx)\nu_1(x)L_1(x)U_1(0)-\right.$$
$$\left.-\int\rho(dx)\int_0^{\infty}\int_0^{\tau}\frac{s}{1!}F_x(ds)\mathbb{V}(x)PW_1(\tau-s)d\tau\right]/\widehat{m},$$
where
$\nu_2(x)=m_2(x)-\mu_3(x)=\frac{2m_1(x)m_2(x)-m_3(x)}{3!m_1(x)}.$

For the next terms we have by analogy:
$$c_k(0)=\left[\sum_{r=0}^{k-1}\int\pi(dx)\nu_{k-r}(x)L_{k-r}(x)U_r(0)-\right.$$ $$\left.-\sum_{r=1}^{k-1}\int\rho(dx)\int_0^{\infty}\int_0^{\tau}\frac{s^r}{r!}F_x(ds)\mathbb{V}^r(x)PW_{k-r}(\tau-s)d\tau\right]/\widehat{m},$$ where $\nu_k(x)=(-1)^{k+1}[\mu_{k+1}(x)-m_k(x)].$

\begin{rem} The last term of this equality may be rewritten without the
use of singular terms. This may be done, like in the work {\rm
\cite{Sam}}, using the Laplace transform of the $W_i(\tau)$.
\end{rem}

This finishes the proof of Theorem 1.2. $\Box$

\section{\bf Estimate of the remainder}

We suppose in addition to the previous assumptions that the function
$v(u;x)$ in (5), (6) satisfies the condition
$$\sup\limits_{|u|\leq R }\sup\limits_{x\in E }|v(u;x)|\leq
C_{R},\eqno(33)$$ where $R>0$.

\begin{lem} If $\Phi^{\varepsilon}$ solves
$$\mathbb{L}^{\varepsilon}\Phi^{\varepsilon}:=[Q+\varepsilon Q_1]\Phi^{\varepsilon}=\varepsilon\psi, \eqno(34)$$
where $\psi$ is given, such that $||\psi||\leq C, ||Q^{-1}||\leq C,
||Q_1||\leq C,$ then
$$||\Phi^{\varepsilon}||\leq C_1\varepsilon.$$
\end{lem}

{\it Proof.} The solution of equation (34) is:
$$\Phi^{\varepsilon}=\varepsilon Q^{-1}[\psi-Q_1\Phi^{\varepsilon}]. $$

We can make the following iteration: set
$$\Phi^{\varepsilon}_{(0)}=\varepsilon\psi,$$ $$||\Phi^{\varepsilon}_{(0)}||=\varepsilon ||\psi||\leq \varepsilon C,$$
then
$$\Phi^{\varepsilon}_{(1)}=\varepsilon Q^{-1}[\psi-Q_1\Phi^{\varepsilon}_{(0)}],$$ as soon as
$Q^{-1}$ and $Q_1$ are bounded
$$||\Phi^{\varepsilon}_{(1)}||\leq\varepsilon C(C+C||\Phi^{\varepsilon}_{(0)}||)\leq C(\varepsilon
C+(\varepsilon C)^2).$$

By induction, we have:
$$\Phi^{\varepsilon}_{(N)}\leq C\sum_{k=1}^N(\varepsilon C)^k.$$

We may choose the parameter $\varepsilon$ small enough, so that
$\varepsilon C<1.$ Thus, for $N\to\infty$ we obtain:
$$||\Phi^{\varepsilon}||\leq C\frac{\varepsilon C}{1-\varepsilon C}\leq \varepsilon C_1.$$

The lemma is thus proved. $\Box$

\begin{cor} If $\Phi^{\varepsilon}$ solves
$$\mathbb{L}^{\varepsilon}\Phi^{\varepsilon}:=[Q+\varepsilon Q_1]\Phi^{\varepsilon}=\varepsilon^{N+1}\psi, \eqno(35)$$
where $||\psi||\leq C, ||Q^{-1}||\leq C, ||Q_1||\leq C,$ then
$$||\Phi^{\varepsilon}||\leq C_1\varepsilon^{N+1}.$$
\end{cor}

{\it Proof.} Let us suppose
$$\Phi^{\varepsilon}:=\varepsilon^N\overline{\Phi}^{\varepsilon}.$$

Then, we have from (35): $$\varepsilon^N
\mathbb{L}^{\varepsilon}\overline{\Phi}^{\varepsilon}=\varepsilon^{N+1}\psi.$$

So, we obtain an equation for $\overline{\Phi}^{\varepsilon}$:
$$\mathbb{L}^{\varepsilon}\overline{\Phi}^{\varepsilon}=\varepsilon\psi.$$

By Lemma 7.1, we have
$$||\overline{\Phi}^{\varepsilon}||\leq C_1\varepsilon,$$ or, for
$\Phi^{\varepsilon}$: $$||\Phi^{\varepsilon}||\leq
C_1\varepsilon^{N+1}.$$

Hence the corollary is proved. $\Box$

Let us consider the estimate of the remainder, written in the form:
$$\Phi^{\varepsilon,N}(t)=\Phi_t^{\varepsilon}(u,x)-\Phi_N^{\varepsilon}(t)=\Phi_t^{\varepsilon}(u,x)-U_0(t)-\sum_{k=1}^N\varepsilon^k(U_k(t)+W_k(\tau)).$$

\begin{thm} Under the conditions (1) and (33) the
following estimate for the remainder in the expansion of Theorem 1.2
holds:
$$||\Phi^{\varepsilon,N}(t)||\leq
C^{\varepsilon}\varepsilon^{N+1},$$ for some $C$.
\end{thm}

 {\it Proof.} To prove this Theorem we should show that
$\Phi^{\varepsilon,N}(t)$ satisfies the conditions of Corollary 7.2.

Let us first consider
$U^{\varepsilon,N}(t):=U^{\varepsilon}(t)-U_0(t)-\sum_{k=1}^N\varepsilon^kU_k(t),$
that is the regular part of $\Phi^{\varepsilon,N}(t)$.

From the equation (13) we obtain:
$$\int_0^{\infty}F_x(ds)\mathbb{V}_{\varepsilon
s}(x)PU^{\varepsilon,N}(t-\varepsilon s)-
U^{\varepsilon,N}(t)=\int_0^{\infty}F_x(ds)\mathbb{V}_{\varepsilon
s}(x)P[U^{\varepsilon,N}(t-\varepsilon s)-U^{\varepsilon,N}(t)]+$$
$$+\int_0^{\infty}F_x(ds)\mathbb{V}_{\varepsilon
s}(x)PU^{\varepsilon,N}(t)-
U^{\varepsilon,N}(t)=QU^{\varepsilon,N}(t)+\int_0^{\infty}F_x(ds)[\mathbb{V}_{\varepsilon
s}(x)-I]PU^{\varepsilon,N}(t-\varepsilon s)-$$
$$-\varepsilon\int_0^{\infty}sF_x(ds)\mathbb{V}_{\varepsilon
s}(x)P(U^{\varepsilon,N}(t))'+\int_0^{\infty}F_x(ds)\mathbb{V}_{\varepsilon
s}(x)P[U^{\varepsilon,N}(t-\varepsilon
s)-U^{\varepsilon,N}(t)+\varepsilon s(U^{\varepsilon,N}(t))']=$$
$$=QU^{\varepsilon,N}(t)-\overline{F}_x(s)[\mathbb{V}_{\varepsilon
s}(x)-I]PU^{\varepsilon,N}(t)|_0^{\infty}+\varepsilon \mathbb{V}(x)
\int_0^{\infty}\overline{F}_x(s)\mathbb{V}_{\varepsilon
s}(x)PU^{\varepsilon,N}(t)ds-$$
$$-\varepsilon\int_0^{\infty}sF_x(ds)\mathbb{V}_{\varepsilon
s}(x)P(U^{\varepsilon,N}(t))'+\int_0^{\infty}F_x(ds)\mathbb{V}_{\varepsilon
s}(x)P[U^{\varepsilon,N}(t-\varepsilon
s)-U^{\varepsilon,N}(t)+\varepsilon s(U^{\varepsilon,N}(t))']=$$
$$=[Q+\varepsilon Q_1^{\varepsilon}]U^{\varepsilon,N}(t),$$ where
$$Q_1^{\varepsilon}U^{\varepsilon,N}(t):=\mathbb{V}(x)
\int_0^{\infty}\overline{F}_x(s)\mathbb{V}_{\varepsilon
s}(x)PU^{\varepsilon,N}(t)ds-\int_0^{\infty}sF_x(ds)\mathbb{V}_{\varepsilon
s}(x)P(U^{\varepsilon,N}(t))'+$$
$$+\int_0^{\infty}F_x(ds)\mathbb{V}_{\varepsilon
s}(x)P[U^{\varepsilon,N}(t-\varepsilon
s)-U^{\varepsilon,N}(t)+\varepsilon
s(U^{\varepsilon,N}(t))']/\varepsilon, \eqno(36)$$ and we have by
conditions (1) and (33):
$$||Q_1^{\varepsilon}||=||\mathbb{V}
\int_0^{\infty}\overline{F}_x(s)\mathbb{V}_{\varepsilon
s}PU^{\varepsilon,N}(t)ds-\int_0^{\infty}sF_x(ds)\mathbb{V}_{\varepsilon
s}(x)P(U^{\varepsilon,N}(t))'+$$
$$+\int_0^{\infty}F_x(ds)\mathbb{V}_{\varepsilon
s}(x)P[U^{\varepsilon,N}(t-\varepsilon
s)-U^{\varepsilon,N}(t)+\varepsilon
s(U^{\varepsilon,N}(t))']/\varepsilon||\leq$$ $$\leq
||\mathbb{V}||||U^{\varepsilon,N}(t)||\int_0^{\infty}e^{\varepsilon
s
||\mathbb{V}||}\overline{F}_x(s)ds+||(U^{\varepsilon,N}(t))'||\int_0^{\infty}sF_x(ds)e^{\varepsilon
s ||\mathbb{V}||}+ $$ $$+\int_0^{\infty}F_x(ds)e^{\varepsilon s
||\mathbb{V}||}||U^{\varepsilon,N}(t-\varepsilon
s)-U^{\varepsilon,N}(t)+\varepsilon
s(U^{\varepsilon,N}(t))'||/\varepsilon\leq$$ $$\leq C_1+C_2+C_3=C,
$$ as soon as the functions $U^{\varepsilon,N}(t)$ are from the
space of real-valued functions which are bounded with all their
derivatives.

The operator $Q$ has bounded inverse operator
$\mathbb{R}_0=[Q+\Pi]^{-1}-\Pi$ in $R_Q$ (see {\rm \cite{KoTu}}).

On the other hand, we have from (13):
$$\int_0^{\infty}F_x(ds)\mathbb{V}_{\varepsilon
s}(x)PU^{\varepsilon,N}(t-\varepsilon s)-
U^{\varepsilon,N}(t)=\left[\int_0^{\infty}F_x(ds)\mathbb{V}_{\varepsilon
s}(x)PU^{\varepsilon}(t-\varepsilon s)- U^{\varepsilon}(t)\right]-$$
$$-
\left[Q+\sum_{k=1}^{N}\varepsilon^k\mu_k(x)L_k\right]\left[U_0+\sum_{k=1}^N\varepsilon^kU_k\right]
-\frac{(\varepsilon
s)^{N+1}}{(N+1)!}\mathbb{V}^{(N+1)}(x)\int_0^{\infty}F_x(ds)
\int_0^{\varepsilon s}\mathbb{V}_{\theta}(x)\times$$
$$\times PU_N^{\varepsilon}(t-\varepsilon s)d \theta=
-\left[QU_0+\sum_{k=1}^N\varepsilon^k\left(QU_k+\sum_{r=1}^k\mu_r(x)L_rU_{N-r}\right)-\varepsilon^{N+1}\psi\right]
=\varepsilon^{N+1}\psi,$$ where the last equality follows from (15),
the function $\psi$ being an integral of the form (36) and may thus
be estimated in the same way as the operator $Q_1^{\varepsilon}$.

So, we showed that the function $U^{\varepsilon,N}(t)$ satisfies the
equation
$$[Q+\varepsilon Q_1^{\varepsilon}]U^{\varepsilon,N}(t)=\varepsilon^{N+1}\psi,\eqno(37)$$ where the conditions of Corollary 7.2
 are satisfied for $Q, Q_1^{\varepsilon},\psi$.

So, we obtain in $R_Q$: $$||U^{\varepsilon,N}(t)||\leq C
\varepsilon^{N+1}.\eqno(38)$$

It is easy to see that equation (37) in $N_Q$ has the form:
$$\varepsilon Q_1^{\varepsilon}c^{\varepsilon,N}(t)=\varepsilon^{N+1}\psi.\eqno(39)$$

From (36) we have:
$$Q_1^{\varepsilon}c^{\varepsilon,N}(t)=m_1(x)[\mathbb{V}(x)
c^{\varepsilon,N}(t)-(c^{\varepsilon,N}(t))']+\varepsilon\mathbb{V}^2(x)
\int_0^{\infty}\overline{F}_x(s)\int_0^{\varepsilon
s}\mathbb{V}_{\theta}(x)Pc^{\varepsilon,N}(t)d\theta ds-$$
$$-\varepsilon\mathbb{V}(x)\int_0^{\infty}sF_x(ds)\int_0^{\varepsilon
s}\mathbb{V}_{\theta}(x)P(c^{\varepsilon,N}(t))'d\theta+\int_0^{\infty}F_x(ds)\mathbb{V}_{\varepsilon
s}(x)P[U^{\varepsilon,N}(t-\varepsilon s)-$$
$$-U^{\varepsilon,N}(t)+\varepsilon
s(U^{\varepsilon,N}(t))']/\varepsilon.$$

We obtain (39) in the following form: $$[Q_1+\varepsilon
Q_2^{\varepsilon}]c^{\varepsilon,N}(t)=\varepsilon^{N}\psi,$$ where
$Q_1^{-1}$ is a bounded integral operator, $Q_2^{\varepsilon}$ and
$\psi$ are integrals of the form (36) and may thus be estimated in
the same way as the operator $Q_1^{\varepsilon}$.

So the conditions of Corollary 7.2 are satisfied for $Q_1,
Q_2^{\varepsilon},\psi,$ and we have in $N_Q$:
$$||c^{\varepsilon,N}(t)||\leq C
\varepsilon^{N}.\eqno(40)$$

Let us now consider
$W^{\varepsilon,N}(\tau):=W^{\varepsilon}(\tau)-\sum_{k=1}^N\varepsilon^kW_k(\tau),$
that is the singular part of $\Phi^{\varepsilon,N}(t)$.

From the equation (13) we have:
$$\int_0^{\infty}F_x(ds)\mathbb{V}_{\varepsilon
s}(x)PW^{\varepsilon,N}(\tau- s)-
W^{\varepsilon,N}(\tau)=\int_0^{\tau}F_x(ds)\mathbb{V}_{\varepsilon
s}(x)PW^{\varepsilon,N}(\tau- s)- W^{\varepsilon,N}(\tau)+ $$ $$+
\int_{\tau}^{\infty}F_x(ds)\mathbb{V}_{\varepsilon
s}(x)PW^{\varepsilon,N}(\tau- s)=
\int_0^{\tau}F_x(ds)PW^{\varepsilon,N}(\tau- s)-
W^{\varepsilon,N}(\tau)+$$ $$+
\int_{\tau}^{\infty}F_x(ds)\mathbb{V}_{\varepsilon
s}(x)PW^{\varepsilon,N}(\tau- s)+
\int_0^{\tau}F_x(ds)[\mathbb{V}_{\varepsilon
s}(x)-I]PW^{\varepsilon,N}(\tau- s)
=[\mathbf{Q}(\tau)-I]W^{\varepsilon,N}(\tau)-$$ $$-
\overline{F}_x(s)\mathbb{V}_{\varepsilon
s}(x)PW^{\varepsilon,N}(\tau- s)|_{\tau}^{\infty} +\varepsilon
\mathbb{V}(x)
\int_{\tau}^{\infty}\overline{F}_x(s)\mathbb{V}_{\varepsilon
s}(x)PW^{\varepsilon,N}(\tau- s)ds-$$ $$-
\overline{F}_x(s)[\mathbb{V}_{\varepsilon
s}(x)-I]PW^{\varepsilon,N}(\tau- s)|_{0}^{\tau}+\varepsilon
\mathbb{V}(x)
\int_{0}^{\tau}\overline{F}_x(s)\mathbb{V}_{\varepsilon
s}(x)PW^{\varepsilon,N}(\tau- s)ds=$$ $$=
[\mathbf{Q}(\tau)-I]W^{\varepsilon,N}(\tau)+\overline{F}_x(\tau)\mathbb{V}_{t}(x)PW^{\varepsilon,N}(0)-\overline{F}_x(\tau)
[\mathbb{V}_{t}(x)-I]PW^{\varepsilon,N}(0)+$$ $$+\varepsilon
\mathbb{V}(x)
\int_{0}^{\infty}\overline{F}_x(s)\mathbb{V}_{\varepsilon
s}(x)PW^{\varepsilon,N}(\tau-
s)ds=[\mathbf{Q}(\tau)-I]W^{\varepsilon,N}(\tau)+\overline{F}_x(\tau)PW^{\varepsilon,N}(0)+$$
$$+ \varepsilon \mathbb{V}(x)
\int_{0}^{\infty}\overline{F}_x(s)\mathbb{V}_{\varepsilon
s}(x)PW^{\varepsilon,N}(\tau-
s)ds=([\mathbf{Q}(\tau)-I]+\varepsilon
\mathbf{Q}_1)W^{\varepsilon,N}(\tau),$$ where due to Lemma 5.1
$$\varepsilon\mathbf{Q}_1W^{\varepsilon,N}(\tau):=\varepsilon\mathbb{V}(x)
\int_0^{\infty}\overline{F}_x(s)\mathbb{V}_{\varepsilon
s}(x)PW^{\varepsilon,N}(\tau-
s)ds+\overline{F}_x(\tau)PU^{\varepsilon,N}(0),
$$ and we have by conditions (1), (33) and (38):
$$||\mathbf{Q}_1||=||\mathbb{V}
\int_0^{\infty}\overline{F}_x(s)\mathbb{V}_{\varepsilon
s}(x)Pds||+||\overline{F}_x(\tau)PU^{\varepsilon,N}(0)||\leq$$
$$\leq ||\mathbb{V}||\int_0^{\infty}e^{\varepsilon s
||\mathbb{V}||}\overline{F}_x(s)ds+||U^{\varepsilon,N}(0)||\leq
C_1+\varepsilon^{N+1}C_2\leq C. $$

The conditions of Chapter 1, sections 3,4 from {\rm \cite{Shu}} are
fulfilled for the operator $\mathbf{Q}(\tau)$. In fact,
$$\int_E\int_{{E}}\int_0^{\infty}\pi(dx)F_x(dt)P(x,dy)=1<\infty.$$
So, the operator $\mathbf{Q}(\tau)-I$ has an inverse operator which
is bounded.

On the other hand, we have:
$$\int_0^{\infty}F_x(ds)\mathbb{V}_{\varepsilon
s}(x)PW^{\varepsilon,N}(\tau- s)-
W^{\varepsilon,N}(\tau)=\left[\int_0^{\infty}F_x(ds)\mathbb{V}_{\varepsilon
s}(x)PW^{\varepsilon}(\tau- s)- W^{\varepsilon}(\tau)\right]-$$ $$-
\left[(\mathbf{Q}-I)+
\sum_{k=1}^{N}\varepsilon^k\int_0^{\infty}F_x(ds)\frac{s^k}{k!}\mathbb{V}^k(x)P\right]
\left[\sum_{k=1}^N\varepsilon^kW_k\right]-$$
$$-\int_0^{\infty}F_x(ds)\frac{\varepsilon^{N+1}\mathbb{V}^{N+1}(x)}{(N+1)!}
\int_0^{\varepsilon
s}\mathbb{V}_{\theta}(x)PW_{N}^{\varepsilon}(\tau- s)d\theta.$$

As soon as the term
$\left[\int_0^{\infty}F_x(ds)\mathbb{V}_{\varepsilon
s}(x)PW^{\varepsilon}(\tau- s)- W^{\varepsilon}(\tau)\right]$
satisfies equation (17), we obtain:
$$\int_0^{\infty}F_x(ds)\mathbb{V}_{\varepsilon
s}(x)PW^{\varepsilon,N}(\tau- s)- W^{\varepsilon,N}(\tau)=
\varepsilon
\mathbb{V}(x)\int_{\tau}^{\infty}\overline{F}_x(s)\mathbb{V}_{\varepsilon
s}(x) \varphi(u)ds - $$ $$-
\sum_{k=1}^{N}\varepsilon^k\left[(\mathbf{Q}-I)W_k(\tau)+\psi_0^k(\tau)\right]+
O(\varepsilon^{N+1})= \sum_{k=1}^N\psi^k(\tau)+O(\varepsilon^{N+1})-
$$ $$-\sum_{k=1}^{N}\varepsilon^k\left[(\mathbf{Q}-I)W_k(\tau)+\psi_0^k(\tau)\right]+O(\varepsilon^{N+1})
=\varepsilon^{N+1}\phi,$$ where we used the equalities (18). The
function $\phi$ is an integral of the form (36) and may be estimated
in the same way as the operator $Q_1$.

So, we showed that the function $W^{\varepsilon,N}(\tau)$ satisfies the
equation
$$([\mathbf{Q}(\tau)-I]+\varepsilon \mathbf{Q}_1)W^{\varepsilon,N}(\tau)=\varepsilon^{N+1}\phi,$$
where $\mathbf{Q}(\tau)-I, \mathbf{Q}_1$ and $\phi$ satisfy the
conditions of Corollary 7.2.

So, we obtain: $$||W^{\varepsilon,N}(\tau)||\leq C
\varepsilon^{N+1}.\eqno(41)$$

For the function $\Phi^{\varepsilon,N}(t)$ we have from (38), (40)
and (41):
$$||\Phi^{\varepsilon,N}(t)||\leq ||U^{\varepsilon,N}(t)||+||c^{\varepsilon,N}(t)||+||W^{\varepsilon,N}(\tau)||\leq C\varepsilon^{N+1}+C\varepsilon^{N}+C\varepsilon^{N+1}=C^{\varepsilon}\varepsilon^{N}.$$

Finally, if we consider
$$\Phi^{\varepsilon,N+1}(t)=\Phi^{\varepsilon,N}(t)-\varepsilon^{N+1}[U^{\varepsilon,N+1}(t)+c^{\varepsilon,N+1}(t)+W^{\varepsilon,N+1}(t)],$$
then
$$\Phi^{\varepsilon,N}(t)=\Phi^{\varepsilon,N+1}(t)+\varepsilon^{N+1}[U^{\varepsilon,N+1}(t)+c^{\varepsilon,N+1}(t)+W^{\varepsilon,N+1}(t)].$$

As a result, the following estimation holds:
$$||\Phi^{\varepsilon,N}(t)||\leq||\Phi^{\varepsilon,N+1}(t)||+O(\varepsilon^{N+1})\leq C^{\varepsilon}\varepsilon^{N+1}.$$

Thus theorem 7.3 is proved. $\Box$

{\it Acknowledgements.}This work was supported by DFG  projects
436 UKR 113/70/0-1 and 436 UKR 113/80/0-1, which is gratefully
acknowledged.\par


\begin{thebibliography}{99}

\bibitem{GH} Griego R., Hersh R., \textit{Random evolutions, markov chains,
and systems of partial differential equations,} Proc. Nat. Acad.
Sci. U.S.A. 62, 305–308 (1969).

\bibitem{H1} Hersh R., \textit{Random evolutions: a survey of results and
problems,} Rocky Mountain J. Math. 4, 443–477 (1974).

\bibitem{H2} Hersh R., \textit{The birth of random evolutions,} Mathematical
Intelligencer \underline{25}(1), 53–60 (2003).

\bibitem{HP} Hersh R., Pinsky M., \textit{Random evolutions are asymptotically Gaussian,}
Comm. Pure Appl. Math. 25, 33-44 (1972).

\bibitem{Hill} Hillen T., \textit{Transport Equations and Chemosensitive Movement,}
Habilitation thesis, University of T\"{u}bingen, T\"{u}bingen,
Germany (2001).

\bibitem{Kor1} Korolyuk V.S., \textit{Stochastic systems with averaging in the scheme of diffusion approximation,}
 Ukrainian Math. Journ. \underline{57}, 9, 1235-1252 (2005).

\bibitem{Kor2} Korolyuk V.S., {\it Boundary layer in asymptotic
analysis for random walks,} Theory of Stochastic Processes
\underline{1-2}, 25-36 (1998).

\bibitem{KoKo} Korolyuk V.S., Korolyuk V.V., \textit{Stochastic Models
of Systems,} Kluwer Acad. Publ. (1999), 250p.

\bibitem{KoLim} Korolyuk V.S., Limnios N., \textit{Stochastic Systems in Merging Phase Space,}
World Scientific Publishers (2005), 330p.

\bibitem{KoPeTu} Korolyuk V.S., Penev I.P., Turbin A.F., {\it Asymptotic
expansion for the distribution of absorption time of Markov chain,}
Cybernetics \underline{4}, 133-135 (1973). (in Russian)

\bibitem{KoTu} Korolyuk V.S., Turbin A.F., \textit{Mathematical foundation of
state lumping of large systems,} Kluwer Acad. Publ. (1990), 280p.

\bibitem{KoTurb} Korolyuk V.S., Turbin A.F., \textit{Semi-Markov processes and applications,}
Naukova dumka, Kiev (1976), 181p.

\bibitem{Kush} Kushner H.J., \textit{Approximation and Weak Convergence Methods for Random
Processes, with Applications to Stochastic Systems Theory,} MIT
Press (1984), 269p.

\bibitem{OH} Othmer H.G., Hillen T., \textit{The diffusion limit of transport equations II:
chemotaxis equations,} Siam J. Appl. Math. \underline{62}(4),
1222–1250 (2002).

\bibitem{P1} Papanicolaou G.C., \textit{Asymptotic analysis of transport processes,} Bull.
Amer. Math. Soc. 81, 330-392 (1975).

\bibitem{P} Papanicolaou G.C., \textit{Probabilistic problems and methods in singular
perturbations,} Rocky Mountain Journal of Math. 6, 653-673 (1976).

\bibitem{Pin} Pinsky M.A., \textit{Lectures on Random Evolutions,} World Scientific (1991),
136p.

\bibitem{Sam}  Samoilenko I.V., \textit{Asymptotic expansion for the
functional of markovian evolution in $R^d$ in the circuit of diffusion
approximation,} Journal of Applied Mathematics and Stochastic Analysis 3,
247-258 (2005).

\bibitem{Shu} Shurenkov V.M., \textit{Ergodic Markov processes,}
Nauka, Moskow (1989), 336p. (in Russian)

\bibitem{Skor} Skorokhod A.V., Hoppensteadt F.C., Salehi H.
\textit{Random Perturbation Methods with Applications in Science
and Engineering,} Springer (2002), 488p.

\bibitem{Tad} Tajiev A., {\it Asymptotic
expansion for the distribution of absorption time of semi-Markov process,}
Ukrainian Math. Journ. 9, 422-426 (1978). (in Russian)

\bibitem{VaBu}  Vasiljeva A.B., Butuzov V.F., \textit{Asymptotic methods in the theory of singular
perturbations,} Vysshaja shkola, Moscow (1990), 208p. (in Russian)

\bibitem{YZ} Yin G.G., Zhang Q., \textit{Continuous-Time Markov Chains and Applications:
a Singular Perturbation Approach,} Springer (1998), 349p.


\end{thebibliography}
\end{document}